\documentclass[11pt
 ]{elsarticle}
\usepackage{amsmath}
\usepackage{amscd}
\usepackage{amsfonts}
\usepackage{amsthm}
\usepackage{graphicx}

\usepackage{color}
\usepackage[colorlinks, linkcolor=red, citecolor=green,filecolor=magenta,urlcolor=cyan]{hyperref}

\def\x{{\cal x}}

\def\ds{\displaystyle}
\newfont{\Blackboard}{msbm10 scaled 1200}

\newfont{\roma}{cmr10 scaled 1200}

\newtheorem{thm}{{}\hskip\parindent Theorem}[section]
\newtheorem{lem}{{}\hskip\parindent Lemma}[section]
\newtheorem{pro}{{}\hskip\parindent Proposition}[section]
\newtheorem{exl}{{}\hskip\parindent Example}[section]

\newtheorem{rem}{{}\hskip\parindent Remark}[section]

\def\be{\begin{equation}}
\def\ee{\end{equation}}
\def\beq{\arraycolsep=1.5pt\begin{eqnarray}}
\def\eeq{\end{eqnarray}}

\journal{...}

\begin{document}
\begin{frontmatter}

\title{Existence of optimal pairs for optimal control problems with states constrained to Riemannian manifolds}

\author[Deng]{Li Deng}
\address[Deng]{School of Mathematics,   Southwest Jiaotong University, Chengdu 611756, China. E-mail
address: dengli@swjtu.edu.cn}

\author[2]{Xu Zhang\corref{cor}}
\address[2]{School of Mathematics and New Cornerstone Science Laboratory,
Sichuan University, Chengdu, 610064, China. E-mail
address: zhang\_xu@scu.edu.cn.}

\cortext[cor]{Corresponding author.
\\  This work is supported by the National Science
Foundation of China under grants 11401491, 11821001, 11931011 and 12371451, by the Fundamental research funds for the Central Universities under grant 2682014CX052, by the New
Cornerstone Science Foundation,  by the Science
Development Project of Sichuan University under grant 2020SCUNL201, and by Natural Science Foundation of Sichuan province under grant 2023NSFSC1342.}

\begin{abstract}
In this paper,  we investigate the existence of optimal pairs for optimal control problems with their states constrained pointwise to Riemannian manifolds. For this purpose, by means of the Riemannian geometric tool, we introduce a crucial Cesari-type property, which is an extension of the classical Cesari property (See \cite[Definition 3.3, p. 51]{Berkovitz}) from the setting of Euclidean spaces to that of Riemannian manifolds. Moreover, we show the efficiency of our result by a concrete example.
\end{abstract}


\end{frontmatter}

\setcounter{equation}{0}

\section{Introduction }
\def\theequation{1.\arabic{equation}}
\hskip\parindent

 Let $n\in\mathbb N$ and $M$ be an $n$ dimensional complete Riemannian manifold with a Riemannian metric $g$. Denote by $T_xM$ the tangent space of $M$ at $x\in M$, by $T M\stackrel{\triangle}{=}\bigcup_{x\in M}T_xM$ the tangent bundle on $M$, and by $\langle\cdot,\cdot\rangle$ and $|\cdot|$ the inner product and the norm over $T_xM$  related to $g$, respectively (See \cite[Chapter 3]{le1} for more material on these notions). Given $T>0$,  a metric space $U$, maps $f:[0,T]\times M\times U\to T M$ and $\Gamma: [0,T]\times M\to 2^U$ (the class of all subsets of $U$), and subsets $Q\subseteq M$ and $S\subseteq M\times M$, we consider the following control system:
\begin{align}\label{se}
\dot y(t)=f(t,y(t),u(t))\quad a.e.\;t\in[0,T],
\end{align}
with the control $u(\cdot)$ and state $y(\cdot)$ satisfying the following constrains:
\begin{align}\label{cons}
\left\{\begin{array}{l}
u(t)\in \Gamma(t,y(t)),\quad a.e.\;t\in[0,T],
\\[2mm] y(t)\in Q,\quad\forall\;t\in(0,T),
\\[2mm] ( y(0),y(T))\in S.
\end{array}\right.
\end{align}
Write $\mathcal U\stackrel{\triangle}{=}\{w:(0,T)\to U\;|\; w(\cdot) \hbox{ is measurable}\}$. A pair $(u(\cdot),y(\cdot))$ is said to be feasible, if $u(\cdot)\in\mathcal U$, $y(\cdot)$ is absolutely continuous, and both (\ref{se}) and (\ref{cons}) are satisfied.

Given maps $f^0:[0,T]\times M\times U\to \mathbb R$ and $h: M\times M\to \mathbb R$,
a pair $(u(\cdot), y(\cdot))$ is said to be  admissible, if   it is feasible and $f^0(\cdot,y(\cdot),u(\cdot))\in L^1(0,T)$.  Hereafter, we denote by $\mathcal P_{ad}$ the set of all admissible pairs, and introduce the cost functional with respect to the control system as follows:
\begin{align*}
J(u(\cdot),y(\cdot))\stackrel{\triangle}{=}\int_0^Tf^0(t,y(t),u(t))dt+h(y(0),y(T)),\quad\forall\;(u(\cdot),y(\cdot))\in\mathcal P_{ad}.
\end{align*}
Our optimal control problem is formulated as follows.
 \begin{description}
\item[(P)]  Find $(\bar u(\cdot),\bar y(\cdot))\in\mathcal P_{ad}$ such that
\begin{align*}
J(\bar u(\cdot),\bar y(\cdot))=\inf_{(u(\cdot),y(\cdot))\in\mathcal P_{ad}}J(u(\cdot),y(\cdot)).
\end{align*}
\end{description}
If the above pair $(\bar u(\cdot),\bar y(\cdot))$ does exist, then we say that $(\bar u(\cdot),\bar y(\cdot))$, $\bar u(\cdot)$ and $\bar y(\cdot)$  are an {\it optimal pair},  {\it optimal control} and {\it optimal trajectory}, respectively. In this paper, we are concerned with the existence of optimal pairs for Problem (P).

The  study of existence of optimal pairs for optimal control problems began in 1960s. L.~Cesari in \cite{cesari66} studied this sort of problems for finite dimensional systems evolved in Euclidean spaces, in which a key condition, now known as the Cesari property, was introduced to guarantee the existence of optimal pairs. Later, for problems with weaker conditions (compared to the classical Cesari property),  there have been a lot of works pursuing further in this direction, see for example \cite{HL69, olech1970, war72, Berkovitz, Gam78, Marc, cesari83, Raymond87,  Raymond90, CO2003, PJ2007, PJ2009} and the references cited therein.

As far as we know, when the states of the control systems are constrained to differentiable manifolds, there are only quite a few existence results on optimal pairs, for which we mention that A.~Agrachev, D.~Barilari and U.~Boscain considered an affine control system (i.e. a system depending linearly on the control) in \cite[Theorem 3.43, p. 89]{abb2020}.

There are two   motivations for us to study the existence of optimal pairs for Problem (P). First,  the study of Problem (P) is of practical importance. As we have explained in \cite[Section 1]{dzgJ}, the states of many optimal control problems are actually constrained to  manifolds.
Next,  the study on the existence of optimal pairs for Problem (P) is far from satisfactory.
 Even though sometimes one may transform Problem (P) to an optimal control problem on some Euclidean space (i.e., embed $M$ isometrically into a Euclidean space by the Nash embedding theorem (see \cite[Theorem 3]{Nash}, then rewrite the control system (\ref{se}) as a controlled ordinary differential equation on that flat space, and
 replace the hidden state constraint $y(t)\in M$ a.e. $t\in[0,T]$ in  (\ref{se}) by a pointwise equality-type constraint $ c(y(t))=0$ a.e. $t\in[0,T]$ for some function $c$), and then apply the known existence results in the literature to the resulting system, the transformed problem may not satisfy  the conditions required in these existence results.   On the other hand, as mentioned before, the results obtained by geometric tools are considerably limited, for which, to the best of our knowledge,  only an affine control system has been studied (See \cite[Theorem 3.43, p. 89]{abb2020} for example). We refer Remark  \ref{roe1} for a more detailed explanation.

The main purpose of this paper is to find some verifiable conditions guaranteeing the existence of optimal pairs for Problem (P). A fundamental difficulty for this aim is caused by the curved state space $M$.  Actually,  a usual way to prove the existence of optimal pairs (in flat spaces) can be described as follows: Choose a minimizing sequence (with respect to the cost functional) in the set of admissible pairs first and then prove the limit point (in a suitable sense) of this sequence is an admissible pair.  For our present problem, the curved state space $M$ makes this procedure hard to proceed.  To overcome this difficulty, we use Riemannian geometry as an analytic tool and prove the existence of optimal pairs under some proper conditions. It is worth mentioning that, the Cesari-type property (A6) (that we shall introduce in Subsection \ref{sub2.1}) is crucial for the existence of optimal pairs, which is an extension of the classical Cesari property (e.g., \cite[Definition 3.3, p. 51]{Berkovitz}) from the setting of Euclidean spaces to that of Riemannian manifolds.  When the manifolds degenerate to Euclidean spaces, our condition (A6) is exactly the classical Cesari condition. Our main existence result is stated in Theorem  \ref{eth}.  Also, we provide a method (in Proposition \ref{vctc}) to verify the above mentioned Cesari-type property, and give an example (See Example \ref{excc}) as an application.

This paper has the following novelties: First, although
Theorem \ref{eth} is an extension of \cite[Theorem 4.2, p. 58]{Berkovitz} from the Euclidean case to the case of differentiable manifolds, in the present work we have not assumed that the pointwise state-constraint set $Q$ is compact (See the condition (A5)). Next,  Proposition  \ref{vctc} looks very similarly to \cite[Proposition 4.3, p. 107]{ly}.  However, the state space discussed in Proposition  \ref{vctc} is curved, while the state space in the latter is a Banach space (which is still flat).  Thirdly, \cite[Theorem 3.43, p.89]{abb2020} exhibits the existence of the minimizing curve connecting two  fixed points on a sub-Riemannian manifold.  When these points are close enough, the problem
is exactly the existence of optimal pairs for an optimal control problem with an affine control system (for which the system is linear with respect to the control variable) on a differentiable manifold and with a quite special cost functional.
 Comparing to \cite[Theorem 3.43, p.89]{abb2020}, our results have two differences:
1)  We study a more general control system, which is allowed to be nonlinear with respect to the control variable;
2) Our control set $U$ is not required to be bounded.
We also mention that, the optimal control problem discussed in   \cite[Theorem 3.43, p. 89]{abb2020}  fulfills the Cesari-type property (A6) (although such a condition was not introduced there), which will be verified by means of Proposition \ref{vctc}.
 In  Remark \ref{roe1}, we shall explain that,  neither \cite[Theorem 4.2, p. 58]{Berkovitz} nor  \cite[Theorem 3.43, p. 89]{abb2020}  works for the optimal control problem in Example \ref{excc}.

The rest of this paper is organized as follows. In Section \ref{mre}, the main results of this paper are stated, and two  examples are also given.  Section \ref{pm} is devoted to a proof of our main results.

\setcounter{equation}{0}
\hskip\parindent \section{Statement of the main results}\label{mre}
\def\theequation{2.\arabic{equation}}

\subsection{Notations and main assumptions}\label{sub2.1}
First,  we introduce some notations. For the Riemannian manifold $(M,g)$, we denote by $\rho(\cdot,\cdot)$, $\nabla$, and $i(x)$ the distance function, Levi-Civita connection, and the injectivity radius at $x\in M$, respectively.   For $x,y\in M$ with $\rho(x,y)<\min\{i(x),i(y)\}$, we denote by $L_{xy}$ the parallel translation of tensors at $x$ from $x$ to $y$, along the unique shortest geodesic connecting $x$ and $y$. For more details of these notions, we refer the readers to \cite{c,cdz}.

For a.e. $t\in (0,T)$ and each $x\in M $, we define a map $\mathcal E:[0,T]\times T M\to 2^{\mathbb R\times T M}$ as follows:
\begin{align*}
\mathcal E(t,x)\stackrel{\triangle}{=}\big\{(x^0,X)\in\mathbb R\times T_x M\;\big|\;\exists\,u\in \Gamma(t,x)\,s.t.\,x^0\geq f^0(t,x,u)\,\textrm{and}\,X=f(t,x,u)\big\},
\end{align*}
and denote by
\begin{align*}
L_{xy}\mathcal E(t,x)\stackrel{\triangle}{=}\big\{(x^0,L_{xy}X)\;\big|\;(x^0,X)\in\mathcal E(t,x)\big\}\quad\textrm{if}\,\rho(x,y)<i(y).
\end{align*}

We introduce the following assumptions that we shall need in the sequel.
\begin{description}
\item[(A1)] $(U,d)$ is a  complete separable metric space;
\item[(A2)] The function $f:[0,T]\times M\times U\to T M$ is measurable  in $(t,x,u)$.  For a.e. $t\in(0,T)$, the map $(x,u)\mapsto f(t,x,u)$ is continuous,  and there exists a positive constant $K>1$, $x_0\in M$ and $\ell(\cdot)\in L^p(0,T)$ with $p>1$ such that
\begin{align}\label{lcf}
&|L_{x_1x_2}f(s,x_1,u)-f(s,x_2,u)|\leq K\rho(x_1,x_2),
\\\label{bf}&|f(s,x_0,u)|\leq \ell(s)
\end{align}
holds for a.e. $s\in (0,T)$ and all $u\in U$ and $x_1,x_2\in M$ with $\rho(x_1,x_2)<\min\{i(x_1),$ $i(x_2)\}$;

\item[(A3)] The function  $f^0:[0,T]\times M\times U\to\mathbb R$ is Borel measurable in $(t,x,u)$ and lower semicontinuous in $(x,u)$, $h:M\times M\to\mathbb R$ is lower semicontinuous, and there exists a positive constant $K$ such that
\begin{align*}
&\inf\{f^0(t,x,u)\;|\; (x,u)\in M\times U, a.e.\; t\in(0,T)\}\geq -K,
\\&\inf\{h(x_1,x_2)\;|\;x_1,x_2\in M\}\geq -K.
\end{align*}
\item[(A4)] The map $\Gamma:[0,T]\times M\to 2^U$  is a set-valued map with closed images. For every $x\in M$, the map $t\mapsto\Gamma(t,x)$ is measurable (See \cite[Definition 8.1.1, p. 307]{AF}), and for a.e. $t\in[0,T]$,  the map $x\mapsto \Gamma(t,x)$ is continuous (See \cite[Definitions 1.4.1, 1.4.2 and 1.4.3,\;pp. 38--40]{AF}), i.e.,  for each $x\in M$ the following two conditions hold:
\begin{description}
\item[(i)] For any neighborhood $\mathcal O$ of $\Gamma(t,x)$, there exists $\delta>0$ such that, for any $y\in M$ with $\rho(x,y)\leq\delta$, it holds that $\Gamma(t,y)\subset \mathcal O$.

\item[(ii)] For any $v\in \Gamma(t,x)$ and for any sequence $\{x_i\}\subset M$ converging to $x$, there exists a sequence $\{v_i\}\subset U$ satisfying $v_i\in \Gamma(t,x_i)$ for each $i\in\mathbb N$ such that, $v_i$ converges to $v$ as $i$ goes to infinity;
\end{description}
\item[(A5)] $Q\subseteq M$ is closed, and $S\subseteq M\times M$ is bounded and closed;
\item[(A6)] For a.e. $t\in(0,T)$, the map $\mathcal E(t,\cdot): M\to 2^{\mathbb R\times T M}$ satisfies the Cesari-type property on $Q$, i.e. the following property
\begin{align}\label{ctpx}
\bigcap_{i(x)>\delta>0}\textrm{cl co}\,\Big(\bigcup_{\rho(y,x)<\delta}L_{yx}\mathcal E(t,y)\Big)=\mathcal E(t,x)
\end{align}
holds for all $x\in Q$,
where $\textrm{cl}\, C$ and $\textrm{co}\, C$ denote the closure and  convex hulls of the set $C$, respectively.
\end{description}

\subsection{Main results}

Our main existence result on optimal pairs for Problem (P) is stated as follows:
\begin{thm}\label{eth}
Suppose that the assumptions  (A1)-(A6)  hold and $\mathcal P_{ad}\not=\emptyset$. Then Problem (P)
admits at least one optimal pair.
\end{thm}

\begin{rem}\label{roe}
Clearly, Theorem \ref{eth} is an extension of \cite[Theorem 4.2, p. 58]{Berkovitz} from the Euclidean case to the case of manifolds.  In particular, the Cesari-type property (A6) is an extension of the  condition in \cite[(3.5), p. 51]{Berkovitz} by means of the ``parallel translation'', which is an important notion in Riemannian geometry. When the state space $M$ degenerates to a Euclidean space, the condition (A6) is exactly  that in \cite[(3.5), p. 51]{Berkovitz}. Moreover, compared to  \cite[Theorem 4.2, p. 58]{Berkovitz}, we do not require the set $Q$ to be compact. Finally,  in the case that the control system (\ref{se}) is linear, the assumption
(A2) holds for many cases, say when the control set $U$ is bounded (which is reasonable for many concrete problems). In fact, in the existing results for linear control systems (except the problems with special cost functional), similar assumption on the control set is always intrinsically needed, e.g.  \cite[Theorem 4.2, p. 58]{Berkovitz} and \cite{cesari83}.
\end{rem}

As that in the flat spaces, the Cesari-type property (A6) is a kind of convexity condition. The following result provides a method to verify such a property.

\begin{pro}\label{vctc} Assume that for a.e. $t\in(0,T)$, one of the following conditions holds:
\begin{description}
\item[(A$3^\prime$)]  $\Gamma(t,\cdot): M\to 2^U$ is upper semicontinuous (See \cite[Definition 2.1, p. 89]{ly}), the map $f(t,\cdot,\cdot): M\times U\to T M$ is uniformly continuous with respect to $u\in U$ and $f^0(t,\cdot,\cdot): M\times U\to\mathbb R$ is uniformly lower semicontinuous in $u\in U$, i.e., for any $y\in M$ and $\epsilon>0$, there exists $\sigma\in(0,i(y))$ such that
\begin{align*}
|L_{y^\prime y}f(t,y^\prime,u^\prime)-f(t,y,u)|<\epsilon
\quad\textrm{and}\quad f^0(t,y^\prime,u^\prime)>f^0(t,y,u)-\epsilon
\end{align*}
hold for any $(y^\prime,u^\prime,u)\in M\times U\times U$ with $\max\{\rho(y,y^\prime),d(u,u^\prime)\}<\sigma$.

\item[(A$4^\prime$)]   $\Gamma(t,x)\equiv\Gamma(t)$ for all $x\in M$,  the maps $f(t,\cdot,\cdot): M\times U\to T M$ and $f^0(t,\cdot,\cdot): M\times U\to\mathbb R$ satisfy the following condition: for any $y\in M$ and $\epsilon>0$, there exists $\sigma\in(0,i(y))$ such that, whenever $(y^\prime,u)\in M\times U$ satisfies $\rho(y,y^\prime)<\sigma$, one has
\begin{align*}
|L_{y^\prime y}f(t,y^\prime,u)-f(t,y,u)|<\epsilon
\quad\textrm{and}\quad f^0(t,y^\prime,u)>f^0(t,y,u)-\epsilon.
\end{align*}

\end{description}
Then, the map $\mathcal E(t,\cdot)$ satisfies the Cesari-type property on $M$, i.e.  (\ref{ctpx}) holds, if and only if  $\mathcal E(t,x)$ is convex and closed for each $x\in M$.
\end{pro}

We shall prove Theorem \ref{eth} and Proposition \ref{vctc} in Section \ref{pm}. The following example is an application of   Theorem \ref{eth} and Proposition \ref{vctc}.

\begin{exl}\label{excc}
For $x=(x_1,x_2,x_3)^\top\in\mathbb R^3$, set
\begin{align*}
f_1(x)=(0,e^{-x_1^2}x_1\sin x_3,x_3)^\top,\quad f_2(x)=(e^{-x_2^2}x_2\sin x_3,0,x_3)^\top.
\end{align*}
Consider the following optimal control problem:
Minimize
\begin{align*}
J(u(\cdot),v(\cdot))=&\int_0^T\left(u^2(t)\left(e^{-2x_1(t)^2}x_1(t)^2\frac{(\sin x_3(t))^2}{x_3(t)^2}+1\right)\right.
\\&\left.+v(t)\sqrt{e^{-2x_2(t)^2}x_2(t)^2\frac{(\sin x_3(t))^2}{x_3(t)^2}+1}\right)dt
\end{align*}
over $(u(\cdot),v(\cdot))\in\mathcal U\stackrel{\triangle}{=}\{(u,v):[0,T]\to[0,\pi]\times[0,1]\;|\; (u,v)(\cdot)\;\textrm{is measurable}\}$, where $x(\cdot)=(x_1(\cdot),x_2(\cdot),x_3(\cdot))^\top:[0,T]\to\mathbb R^3$ is the  solution to the following control system:
\begin{align}\label{ecs}
\left\{\begin{array}{l}\dot x(t)=f_1(x(t))\sin u(t)+ f_2(x(t))v(t),\quad a.e.\;t\in[0,T],
\\[3mm] x_3(t)>0,\quad\forall\;t\in[0,T],
\\[3mm] x(0)=(0,0,1)^\top,\quad e^T\leq  x_3(T)\leq 2e^T.
\end{array}
\right.
\end{align}
Then, as we shall show below, Problem (P) for the above optimal control system admits at least one optimal pair.

Indeed, we consider  the set   $$H^3=\big\{(x_1,x_2,x_3)^\top\in\mathbb R^3\;\big|\;x_3>0\big\}$$ as a manifold, with only one coordinate chart $(H^3, I|_{H^3})$, where $I: \mathbb R^3\to \mathbb R^3$ is the identity map and $I|_{H^3}$ is the restriction of $I$ to the subset $H^3$. We also introduce on  $H^3$ the metric $g$:
\begin{align*}
g_{ij}(x_1,x_2,x_3)=\delta_{ij}/x_3^2,\quad i,j=1,2,3,
\end{align*}
where $\delta_{ii}=1$ for $i=1,2,3$, and $\delta_{ij}=0$ if $i\not=j$.
 The Riemannian manifold $(H^3, g)$  is called the  hyperbolic space of dimension $3$. It was proved in \cite[p. 162]{c} that $(H^3, g)$ is a simply connected and complete Riemannian manifold.

It is obvious that, for each $x\in H^3$, $f_1(x),f_2(x)\in T_x H^3$, and their norms
 in the metric $g$ are as follows:
\begin{align*}
|f_i(x)|_g=\left(\frac{e^{-2x_i^2}x_i^2(\sin x_3)^2+x_3^2}{x_3^2}\right)^{1/2},\quad i=1,2,
\end{align*}
where $|\cdot|_g$ denotes the norm of a tensor with respect to the metric $g$.
Thus, we can view the control system (\ref{ecs}) as a system on  manifold $H^3$:
\begin{align*}
\left\{\begin{array}{l}
\dot{x}(t)=f_1(x(t))\sin u(t) + f_2(x(t))v(t), \quad \text { a.e. } t \in[0, T], \\[3mm]
x(t)\in H^3, \quad \forall\; t \in[0, T], \\[3mm]
x(0)=(0,0,1)^{\top}, \quad e^T\leq x_3(T) \leq 2 e^T,
\end{array}\right.
\end{align*}
with the control $(u(\cdot),v(\cdot))\in\mathcal U$. The corresponding  cost functional can be rewritten as follows:
\begin{align*}
J(u(\cdot),v(\cdot))=\int_0^T\big(u^2(t)|f_1(x(t))|_g^2+v(t)|f_2(x(t))|_g\big)dt.
\end{align*}
In what follows, we will verify  that  this optimal control problem fulfills the assumptions required in Theorem \ref{eth}.

It is easy to check that  (A1), (A3), (A4), (A5), $\mathcal P_{ad}\not=\emptyset$ and (\ref{bf}) are satisfied. It suffices to verify  (\ref{lcf}) and (A6).

First, we shall use \cite[Lemma 4.1]{cdz} to  prove that (\ref{lcf}) holds. The local expression for the corresponding Levi-Civita  connection $\nabla$ is given by the  Christoffel symbols $\Gamma_{ij}^k$, i.e. $\nabla_{\frac{\partial}{\partial x_i}}\frac{\partial}{\partial x_j}=\sum_{k=1}^3\Gamma_{ij}^k\frac{\partial}{\partial x_k}$ for $i,j=1,2,3$. By some direct computations  we obtain  the following (See \cite[p. 161]{c}):
\begin{align}\label{ecs1}
\Gamma_{ij}^k(x)= \frac{1}{x_3}\left(-\delta_{j k} \delta_{i3}-\delta_{k i} \delta_{j3}+\delta_{i j} \delta_{k3}\right), \quad i,j,k=1,2,3,\;x\in H^3.
\end{align}

Assume
$\nabla f_1(x)=\sum_{i,j=1}^3A_{ij}(x)\frac{\partial}{\partial x_i}\otimes dx_j$ for $x\in H^3$. Then, by the definitions of covariant derivative of tensors (See \cite[p.57]{pth}) and Levi-Civita connection (See \cite[Theorem 2.2.2, p.53]{pth}), we have
\begin{align*}
A_{ij}(x)&=\nabla f_1\left(dx_i,\frac{\partial}{\partial x_j}\right)
\\&=\nabla_{\frac{\partial}{\partial x_j}}\left(e^{-x_1^2}x_1\sin x_3\frac{\partial}{\partial x_2}+x_3\frac{\partial}{\partial x_3}\right)(dx_i)
\\&=\frac{\partial}{\partial x_j}\left(e^{-x_1^2} x_1 \sin x_3\right) \delta_{2i}+e^{-x_1^2} x_1 \sin  x_3  \Gamma_{j 2}^i(x)
 +\delta_{j3} \delta_{3i}+x_3 \Gamma_{j 3}^i(x).
\end{align*}
Before estimating the norm of $|\nabla f_1|_g$, we are going to calculate the norms of a vector and covector at $x\in H^3$. Taking any $X=\sum_{i=1}^3X_i\frac{\partial}{\partial x_i}\in T_x H^3$, we have $|X|_g=\frac{1}{x_3}|X|_{\mathbb R^3}$, where $|\cdot|_{\mathbb R^3}$ denotes the Euclidean norm of $X$. Take any covector $\eta=\sum_{i=1}^3\eta_i dx_i$.
According to the definition for the norm of a tensor (See \cite[(2.5)]{cdz}),  we have
\begin{align*}
|\eta|_g&=\sup\{\eta(Y)\;|\;Y=\sum_{i=1}^3Y_i\frac{\partial}{\partial x_i}\in T_x H^3, |Y|_g\leq 1\}
\\&=\sup\left\{\sum_{i=1}^3\eta_i Y_i\;\left|\;\frac{1}{x_3}|(Y_1,Y_2,Y_3)^\top|_{\mathbb R^3}\leq 1\right\}\right.
\\&= x_3|(\eta_1,\eta_2,\eta_3)^\top|_{\mathbb R^3}.
\end{align*}
By means of  \cite[(2.5)]{cdz} again, we also have
\begin{align*}
|\nabla f_1|_g&=\sup\big\{\nabla f_1(\eta,Y)\;\big|\;|\eta|_g\leq 1,|Y|_g\leq 1\big\}
\\&=\sup\big\{\sum_{i,j=1}^3A_{ij}(x)\eta_iX_j\;\big|\;|(\eta_1,\eta_2,\eta_3)^\top|_{\mathbb R^3}\leq 1/x_3,\;|(X_1,X_2,X_3)^\top|_{\mathbb R^3}\leq x_3\big\}
\\&\leq C_1,
\end{align*}
for a big constant $C_1$,
where we have used (\ref{ecs1})
 and the expressions for $A_{ij}(x)$. Similarly, we can also show $|\nabla f_2(x)|_g\leq C_1$ for any $x\in H^3$.
 Thus, (\ref{lcf}) follows from \cite[Lemma 4.1]{cdz}.

Next, we will show that (A6) holds by means of  Proposition  \ref{vctc}. It is obvious that the assumptions $(A3^\prime)$ and $(A4^\prime)$ hold. We claim that, for a.e.  $t\in(0,T)$ the set
 \begin{align*}
\mathcal{E}(t, x)
=\bigcup_{\substack{(u, v) \\
\in[0, \pi] \times[0,1]}}\left\{\left.\left(\begin{array}{c}
z^0 \\
\sin u f_1(x)+v f_2(x)
\end{array}\right)\; \right| \;\begin{array}{l}
(u, v) \in[0, \pi] \times[0,1],
\\
z^0 \geq u^2|f_1(x)|_g^2+v\left|f_2(x)\right|_g
\end{array}
\right\}
\end{align*}
 is convex and closed for each $x\in H^3$.

 It is obvious that $\mathcal E(t,x)$ is closed. We only have to show that it is convex. To this end, take  any
 \begin{align*}
\left(\begin{array}{c}
z_0 \\
\sin u f_1(x)+v f_2(x)
\end{array}\right),\left(\begin{array}{c}
\xi_0 \\
\sin \hat{u} f_1(x)+\hat{v} f_2(x)
\end{array}\right) \in \mathcal{E}\left(t, x\right)
\end{align*}
and $\lambda\in[0,1]$, where
 $(u,v), (\hat u,\hat v)\in[0,\pi]\times [0,1]$. We set $u_\lambda=\lambda u+(1-\lambda)\hat u$ and $v_\lambda=\lambda v+(1-\lambda)\hat v$. Since $\sin x$ with $x\in [0,\pi]$ is concave, we have
 \begin{align*}
\lambda \sin u+(1-\lambda) \sin\hat u\leq \sin u_\lambda.
\end{align*}
On the other hand,  there exists $\hat u_\lambda\in[\min\{u,\hat u\},\max\{u,\hat u\}]$ such that
\begin{align*}
\lambda \sin u+(1-\lambda) \sin\hat u=\sin\hat u_\lambda.
\end{align*}
Let $\hat u_{\lambda,1},u_{\lambda, 1}\in[0,\pi/2]$ be such that
 $\sin \hat u_{\lambda,1}=\sin\hat u_\lambda$ and $\sin u_{\lambda,1}=\sin u_\lambda$. Then $0\leq \hat u_{\lambda,1}\leq u_{\lambda,1}\leq \pi$. Moreover, by the convexity of the function $s^2$ with $s\geq 0$, we further have
 \begin{align*}
\hat{u}_{\lambda,1}^2 \leq u_\lambda^2 \leq \lambda u^2+(1-\lambda) \hat{u}^2.
\end{align*}
Consequently, we obtain
\begin{align*}
& \lambda\left(u^2\left|f_1(x)\right|_g^2+v|f_2(x)|_g\right)+(1-\lambda)\left(\hat{u}\left|f_1(x)\right|_g^2+\hat{v}\left|f_2(x)\right|_g\right) \\
& \geq \hat{u}_{\lambda, 1}^2\left|f_1(x)\right|_g^2+v_\lambda\left|f_2(x)\right|_g,
\end{align*}
with
\begin{align*}
 \lambda\left(\sin u f_1(x)+v f_2(x)\right)+(1-\lambda)\left(\sin \hat u f_1(x)+ \hat{v} f_2(x)\right)  =\sin \hat{u}_{\lambda, 1} f_1(x)+v_\lambda f_2(x).
\end{align*}
Hence, we derive
\begin{align*}
\lambda\left(\begin{array}{c}
z_0 \\
\sin u f_1(x)+v f_2(x)
\end{array}\right)+(1-\lambda)\left(\begin{array}{c}
\xi_0 \\
\sin \hat{u} f_1(x)+\hat{v} f_2(x)
\end{array}\right) \in \mathcal{E}\left(t, x\right),
\end{align*}
which implies $\mathcal E(t,x)$ is convex. 
\end{exl}

\begin{rem}\label{roe1} There are two possible ways to solve Example \ref{excc}. The first one is to use the existing results  for the setting of Euclidean spaces. \cite[Theorem 4.2, p. 58]{Berkovitz} is about   the problem where  the state is constrained to  a compact set at any time.  However,  this result is not applicable to Example \ref{excc}, because the pointwise state constraint set $H^3$ is not compact. Another way is to eliminate the pointwise state constraints in form by transforming the problem to the one evolved on a manifold and using the results for the case of manifolds. In Example \ref{excc} we adopt the latter one.  For the existence of optimal pairs for optimal control problems on differentiable manifolds,   \cite[Theorem 3.43, p.89]{abb2020} concerns on an affine control system with a special cost functional. However, the optimal control problem in Example \ref{excc} does not fulfill
  the assumptions required in  \cite[Theorem 3.43, p.89]{abb2020}.
\end{rem}

In Example  \ref{excc}, the state space $H^3$ has a unique coordinate chart. In what follows, we exhibit another example, in which the state space contains more than one coordinate chart.

\begin{exl}\label{ses}
For  $S^2\stackrel{\triangle}{=}\big\{(x_1,x_2,x_3)^\top\in\mathbb{R}^3\;\big|\;x_1^2+x_2^2+x_3^2= 1\big\}$,  $N= ( 0, 0, 1) ^\top$  and $C>0$ being big enough, we consider the following control system

\begin{align}\label{csse}
\begin{cases}\displaystyle\dot{x}(t)=\sum_{i=1}^3u_i(t)f_i(x(t)),\quad a.e.\,t\,\in(0,T),\\x(0)=N,\quad|x(T)-N|_{\mathbb R^3}\geq \sqrt{2},
\\[3mm]\displaystyle
u(t)=(u_1(t),u_2(t),u_3(t))^\top\in\mathbb R^3,\quad |u(t)|_{\mathbb R^3}\leq C,\quad a.e.\,t\,\in(0,T),
\\[3mm]\displaystyle x(t)\in S^2,\quad\forall\; t\;\in[0,T],&\end{cases}
\end{align}
where the control function $u(\cdot)=(u_1(\cdot),u_2(\cdot),u_3(\cdot))^\top$ is required to be measurable, and $f_1,f_2,f_3$ are maps from $\mathbb{R}^3$ to $\mathbb{R}^3$ given by
\begin{align*}
f_{1}(x) =(x_2,-x_1,0)^\top, \quad
f_{2}(x) =(0,x_3,-x_2)^\top, \quad
f_3(x) =(-x_3,0,x_1)^\top,
\end{align*}
for all $x= ( x_1, x_2, x_3) ^\top\in \mathbb{R} ^3.$ Let $f^0( = f^0( t, x, u) ) : [ 0, T] \times \mathbb{R}^3\times \mathbb R^3\to \mathbb{R} $   be a map satisfying: $1)\,f^0:[0,T]\times \mathbb{R}^3\times \mathbb R^3\to\mathbb{R}$ is measurable in $t\in[0,T]$ and of class $C^1\operatorname{in}\left(x,u\right)\in \mathbb{R}^3\times \mathbb R^3$;\,2) There exists a positive constant $K$ such that $\inf\{f^0(t,x,u)\;|\; x\in S^2,\,|u|_{\mathbb R^3}\leq C,\,a.e.t\in(0,T)\}\geq -K$; 3) For a.e. $t\in(0,T)$ and each $x\in S^2$, the map $f^0(t,x,\cdot)$ is convex. Then, we shall show below the existence of optimal pairs for the following optimal control problem: Minimize $J(u(\cdot))\stackrel{\triangle}{=}\int_0^Tf^0(t,x(t),u(t))dt$ over the pairs $(u(\cdot),x(\cdot))$ satisfying (\ref{csse}).



First,  the control system (\ref{csse}) is  a system evolved on the submanifold $S^2$ of $\mathbb R^3$. In fact, it is obvious that, for each $x\in S^2$, it holds that $x$ is orthogonal to $f_i(x)$ (with $i=1,2,3$). Thus $\{f_1(x),f_2(x),f_3(x)\}\subset T_x S^2$.

Second, we show the existence of admissible pairs. By the same argument as that will be used in ``Step 4'' of the ``Proof of Theorem 2.1'', we only have to show that, for each $x\in S^2$, $\textrm{span}\,\{f_1(x),f_2(x),f_3(x)\}=T_xS^2$, where $\textrm{span}\,\{f_1(x),f_2(x),f_3(x)\}$ is  the set of all linear combinations of $f_1(x),f_2(x)$ and $f_3(x)$.  We will prove this via coordinate charts. The manifold $S^2$ has a family of
  coordinate charts: $(O_+,\varphi_+)$ and $(O_-,\varphi_-)$, where $O_+=S^2\setminus\{(0,0,-1)^\top\}$, $O_-=S^2\setminus\{N\}$ and maps $\varphi_+=(\xi_1,\xi_2)^\top:O_+\to\mathbb R^2$ and $\varphi_-=(\eta_1,\eta_2)^\top: O_-\to\mathbb R^2$ are respectively defined by
\begin{equation}\label{elc}\begin{array}{ll}\displaystyle\xi(x)&\displaystyle=(\xi_1(x),\xi_2(x))^\top=\varphi_+(x)=\left(\frac{x_1}{1+x_3},\frac{x_2}{1+x_3}\right)^\top,\quad\forall\, x\in O_+,
\\[2mm]\displaystyle\eta(x)&\displaystyle=(\eta_1(x),\eta_2(x))^\top=\varphi_-(x)=\left(\frac{x_1}{1-x_3},\frac{x_2}{1-x_3}\right)^\top,\quad\forall\, x\in O_-.
\end{array}\end{equation}
The inverse of $\varphi_+$ and $\varphi_-$ are as follows:
\begin{equation}\label{icc}\begin{array}{ll}
\displaystyle\varphi_+^{-1}({\xi})&\displaystyle=(x_1( {\xi}),x_2( {\xi}),x_3( {\xi}))^\top=\left(\frac{2 {\xi}_1}{a( {\xi})},\frac{2 {\xi}_2}{a( {\xi})},\frac{1-{\xi}_1^2- {\xi}_2^2}{a( {\xi})}\right)^\top,
\\[2mm]\displaystyle {\varphi}_-^{-1}( {\eta})&\displaystyle=(x_1( {\eta}),x_2( {\eta}),x_3( {\eta}))^\top=\left(\frac{2 {\eta}_1}{a( {\eta})},\frac{2 {\eta}_2}{a( {\eta})},\frac{ {\eta}_1^2+ {\eta}_2^2-1}{a( {\eta})}\right)^\top,\end{array}\end{equation}
  for all $\xi=(\xi_1,\xi_2)^\top,\eta=(\eta_1,{\eta}_2)^\top\in\mathbb{R}^2$, where $a:\mathbb{R}^2\to\mathbb{R}$ is given by $a(\xi)=1+ {\xi}_1^2+ {\xi}_2^2$.

For $x\in O_+$, $\left.\{\left.\frac\partial{\partial\xi_1}\right|_x,\left.\frac\partial{\partial\xi_2}\right|_x\}\right.$ is the basis for $T_xS^2$. By means of (\ref{elc}) and (\ref{icc}), we obtain the expressions for $f_1,f_2$ and $f_3$ in the chart $(O_+,\varphi_+)$:
\begin{align*}
&f_{1}( {\xi}) =\xi_2\frac\partial{\partial\xi_1}-\xi_1\frac\partial{\partial\xi_2},  \\
&f_{2}( {\xi}) =\xi_1\xi_2\frac\partial{\partial\xi_1}+\frac{1-\xi_1^2+\xi_2^2}2\frac\partial{\partial\xi_2},  \\
&f_3( {\xi}) =\frac{-1- {\xi}_1^2+ {\xi}_2^2}2\frac\partial{\partial {\xi}_1}- {\xi}_1 {\xi}_2\frac\partial{\partial {\xi}_2},
\end{align*}
which immediately implies $\textrm{span}\{f_1(x),f_2(x),f_3(x)\}=T_xS^2$. When $x\in O_-$, we can obtain this relation similarly.

Third, we show that the condition (A6) holds by means of Proposition  \ref{vctc}.  For a.e. $t\in(0,T)$ and each $x\in S^2$, set
\begin{align*}
\mathcal E(t,x)=\{&(\zeta, X)\in\mathbb R\times T_xS^2\;|\:\exists\, u=(u_1,u_2,u_3)^\top\in\mathbb R^3\,s.t.\,|u|_{\mathbb R^3}\leq C
\\&\textrm{and}\;X=\sum_{i=1}^3u_if_i(x),\;\zeta\geq f^0(t,x,u)\}.
\end{align*}
We obtain from the convexity of $f^0(t,x,\cdot)$ that $\mathcal E(t,x)$ is convex. It is obvious that $\mathcal E(t,x)$ is closed.
 Thus, according to Proposition  \ref{vctc}, we derive that the condition (A6) is satisfied.

For the present optimal control problem, the conditions (A1), (A3), (A4) and (A5) hold obviously. Finally, we show that the condition (A2) holds. We view $S^2$ as a submanifold of $\mathbb R^3$ with induced metric $g$. Since $f_1,f_2,f_3$ are smooth, $S^2$ is compact and the control set is bounded, there exists a positive constant $C_1$ such that
$
|\nabla\sum_{i=1}^3u_if_i(x)|_g\leq C_1$ holds for all $x\in S^2$ and $(u_1,u_2,u_3)^\top\in\mathbb R^3$ with $ |(u_1,u_2,u_3)^\top|_{\mathbb R^3}\leq C$ (Recall that $\nabla$ and $|\cdot|_g$ stand for respectively the Levi-Civita connection and the norm of a tensor with respect to the metric $g$).
Thus, according to \cite[Lemma 4.1]{cdz}, we obtain that, for the present problem (\ref{lcf}) holds. Relation (\ref{bf}) also holds for this problem, since $S^2$ is compact
 and the control set is bounded. 
 \end{exl}

 \begin{rem}\label{roe112}
 We point out that, neither
 \cite[Theorem 4.2, p. 58]{Berkovitz} nor \cite[Theorem 3.43, p. 89]{abb2020} can be used directly to show the existence result in Example  \ref{ses}. Indeed, for the optimal control problem discussed in Example  \ref{ses},  even though the pointwise state constraint set $S^2$ is compact, the existence of the corresponding minimizing sequence of admissible pairs (with an additional property of
equi-absolutely continuous trajectories, required  as an assumption in \cite[Theorem 4.2, p. 58]{Berkovitz}) is unclear. On the other hand,   in \cite[Theorem 3.43, p. 89]{abb2020},  the  state of the  control system is constrained to a  fixed point at the  terminal time, while in control system (\ref{csse}), the state at the terminal time is constrained to a subset of $S^2$.
 \end{rem}

\setcounter{equation}{0}
\hskip\parindent \section{Proof of the main results}\label{pm}
\def\theequation{3.\arabic{equation}}

This section is addressing to proving Theorem \ref{eth} and Proposition \ref{vctc}.

 First, write
 $$C([0,T];M)\stackrel{\triangle}{=}\big\{y:[0,T]\to M\;\big|\; y(\cdot)\,\textrm{is continuous in }M\big\}.$$
We show below some properties of trajectories of the control system  (\ref{se}).

\begin{lem}\label{wpo}
Assume the conditions (A1) and (A2) hold. Then,   for any $x\in  M$ and $u(\cdot)\in \mathcal U$, the system  (\ref{se}) admits a unique solution $y(\cdot)\in C([0,T];M)$ with $y(0)=x$ and
\begin{equation}\label{csd1}
\begin{array}{ll}
\displaystyle\rho(y(s),y(t))\leq \rho(y(t),x_0)\big(e^{K(s-t)}-1\big)+K \int_t^s\ell(\tau)d\tau e^{K(s-t)},\\[3mm]
\displaystyle\qquad\qquad\qquad\forall\;t,s\in [0,T]\hbox{ with }0\leq t<s\leq T,
\end{array}
\end{equation}
where the positive constant $K$ is given in ($A2$).
\end{lem}

{\it Proof.}\; First, we claim that
\begin{align}\label{ef}
|f(t,y,v)|\leq K\rho(x_0,y)+\ell(t),\quad\forall\, (y,v)\in M\times U\;\textrm{and}\;a.e.\,t\in[0,T].
\end{align}
In fact, since $(M,g)$ is complete, it follows from \cite[Theorem 2.8, p. 146]{c} that, for any  $y\in M$, there exists a geodesic $\gamma:[0,1]\to M$ satisfying $\gamma(0)=x_0$, $\gamma(1)=y$, and the length of $\gamma(\cdot)$ (i.e., $\int_0^1|\dot\gamma(t)|dt$) equals to $\rho(x_0,y)$. Moreover, there exist $N\in\mathbb N$ and $\tau_0,\tau_1,\cdots,\tau_N\in [0,1]$ with $0=\tau_0<\tau_1<\cdots<\tau_{N-1}<\tau_N=1$ such that $\rho(\gamma(\tau_j),\gamma(\tau_{j+1}))<\min\{i(\gamma(\tau_j)),i(\gamma(\tau_{j+1}))\}$ for $j=0,1,\cdots,N-1$. For a.e. $t\in[0,T]$ and every $v\in U$, we obtain from \cite[(2.6)]{cdz} (or \cite[ Definition 2.5, p. 52 \& Proposition 3.2, p. 53]{c}), the triangle inequality of the norm $|\cdot|$  and the condition (A2) that
\begin{align*}
|f(t,y,v)|\leq &|f(t,\gamma(\tau_N),v)-L_{\gamma(\tau_{N-1})\gamma(\tau_N)}f(t,\gamma(\tau_{N-1}),v)|+|f(t,\gamma(\tau_{N-1}),v)|
\\ \leq & |f(t,\gamma(\tau_N),v)-L_{\gamma(\tau_{N-1})\gamma(\tau_N)}f(t,\gamma(\tau_{N-1}),v)|
\\&+|f(t,\gamma(\tau_{N-1}),v)-L_{\gamma(\tau_{N-2})\gamma(\tau_{N-1})}f(t,\gamma(\tau_{N-2}),v)|+|f(t,\gamma(\tau_{N-2}),v)|
\\ \leq &\cdots
\\ \leq & \sum_{j=1}^N|f(t,\gamma(\tau_j),v)-L_{\gamma(\tau_{j-1})\gamma(\tau_j)}f(t,\gamma(\tau_{j-1}),v)|+|f(t,x_0,v)|
\\ \leq & K \sum_{j=1}^N \rho(\gamma(\tau_j),\gamma(\tau_{j-1}))+\ell(t)
\\=&K\rho(x_0,y)+\ell(t),
\end{align*}
which implies  (\ref{ef}).

Next, take any $x\in M$ and $u(\cdot)\in\mathcal U$. By employing a local coordinate transformation near $x$,  we can transform the system (\ref{se}) with initial state $y(0)=x$ into an ordinary differential equation with its state valued in $\mathbb R^n$. Applying \cite[Proposition 5.3, p. 66]{ly} to this equation, we obtain that this equation admits a local solution. We now use the contradiction argument to show that such a solution can be extended globally. For this, assume that the system (\ref{se}) admited a solution $y(\cdot)$ on $[0,\tau)$  for some $\tau>0$ and
\begin{align}\label{csd2}
\lim_{s\to\tau^-}\rho(y(s),x)=+\infty.
\end{align}
Take any $t\in[0,\tau)$. For any $s\in[t,\tau)$, we could obtain from  (\ref{ef})  and \cite[Proposition 2.5, p. 146]{c} that
\begin{align*}
&\int_t^s|f(\zeta,y(\zeta),u(\zeta))|d\zeta
\\  & \leq K\int_t^s\big(\rho(y(\zeta),x_0)+\ell(\zeta)\big)d\zeta
\\ &\leq K\int_t^s\rho(y(\zeta), y(t))d\zeta+K\rho(x_0,y(t))(s-t)+K\int_t^s\ell(\zeta)d\zeta.
\end{align*}
By the definition of distance functions on Riemannian manifolds (See \cite[Definition 2.4, p. 146]{c}), it follows that
\begin{align*}
\rho(y(s),y(t))\leq \int_t^s|\dot y(\zeta)|d\zeta=\int_t^s|f(\zeta,y(\zeta),u(\zeta))|d\zeta.
\end{align*}
Applying Gronwall's inequality to the above, we obtain that the inequality (\ref{csd1}) holds  for $s\in[t,\tau)$. To conclude the proof, it suffices to show  $T<\tau$. In fact, if this was not true, we would derive that
\begin{align*}
\lim_{s\to\tau^-}\rho(y(s),y(t))\leq \rho(y(t),x_0)(e^{K(\tau-t)}-1)+K\int_t^\tau \ell(\zeta)d\zeta e^{K(\tau-t)}<+\infty,
\end{align*}
which contradicts (\ref{csd2}).
$\Box$

 \medskip

 \textbf{Proof of Theorem \ref{eth}} The proof is divided into five steps.

\medskip

 {\it Step 1.}\; Since $\mathcal P_{ad}\not=\emptyset$, there exists a sequence $\{(u_k(\cdot),y_k(\cdot))\}\subset\mathcal P_{ad}$ such that $\lim_{k\to+\infty}J(u_k(\cdot),y_k(\cdot))=\inf_{(u(\cdot),y(\cdot))\in\mathcal P_{ad}}J(u(\cdot),y(\cdot))$.  We obtain from the assumption (A5), the completeness of $M$, and  \cite[Theorem 2.8, p. 146]{c} that, there exists a subsequence of $\{y_k(\cdot)\}$ (we still denote it by $\{y_k(\cdot)\}$) and $(y_0,y_T)\in S$ such that $\lim_{k\to+\infty}(y_k(0),y_k(T))$ $=(y_0,y_T)$.
Moreover, it  follows from Lemma \ref{wpo}  and \cite[Proposition 2.5, p. 146]{c} that,
\begin{align*}
\rho(y_k(t),y_k(0))&\leq \rho(y_k(0),x_0)(e^{Kt}-1)+K\int_0^t\ell(\tau) d\tau e^{Kt}
\\&\leq\Big( \rho(y_k(0),y_0)+\rho(y_0,x_0)\Big)(e^{KT}-1)+K\int_0^T\ell(\tau) d\tau e^{KT}
\end{align*}
holds for all $k\in\mathbb N$ and $t\in[0,T]$, which together with the convergence of $\{y_k(0)\}$, immediately implies $\{y_k(\cdot)\}$ is bounded in $C([0,T]; M)$.
    Applying Lemma \ref{wpo} to $\{y_k(\cdot)\}$ again, we obtain that, for all $s\in[0,T)$, it holds that
\begin{align*}
\rho(y_k(t),y_k(s))\leq \rho(y_k(s), x_0)\big(e^{K(t-s)}-1\big)+K \int_s^t\ell(\tau)d\tau e^{K(t-s)}\quad\forall\,t\in[s,T].
\end{align*}
Thus $\{y_k(\cdot)\}$ is equicontinuous in $C([0,T];M)$. Applying Ascoli's theorem(ref. \cite[Theorem 47.1, P.  290]{munk}) to $\{y_k(\cdot)\}$, there exists a subsequence of $\{y_k(\cdot)\}$ (we still denote it by $\{y_k(\cdot)\}$) and $\bar y(\cdot)\in C([0,T];M)$ such that  $y_k(\cdot)$ tends to $\bar y(\cdot)$ in $C([0,T];M)$ as $k\to\infty$.  Recalling that $\{(y_k(0),y_k(T))\}$ converges to $(\bar y(0),\bar y(T))$,  we have $(\bar y(0),\bar y(T))=(y_0,y_T)$. Moreover, since $M$ is complete, we obtain from \cite[Theorem 2.8, p. 146]{c} that, there exists $\bar \delta>0$ such that
\begin{align*}
\inf\big\{i(\bar y(t))\;\big|\;t\in[0,T]\big\}>\bar \delta.
\end{align*}
Then, there exists $\bar k>0$ such that
\begin{align*}
\rho(y_k(t),\bar y(t))\leq \bar \delta<i(\bar y(t)),\quad\forall\,t\in[0,T]\;\textrm{and}\;k\geq\bar k.
\end{align*}
Thus, according to the definition of the parallel translation (e.g., \cite[Section 2.2]{cdz}), we can define
\begin{align*}
L_{y_k(t)\bar y(t)} f(t,y_k(t),u_k(t))\in T_{\bar y(t)} M\quad \textrm{for all} \;k\geq \bar k\;\textrm{and a.e.}\; t\in[0,T].
\end{align*}

\medskip

{\it Step 2.}\;
For $k\geq\bar k$, we express $L_{y_k(\cdot)\bar y(\cdot)}f(\cdot,y_k(\cdot),u_k(\cdot))$ by a frame along $\bar y(\cdot)$. For this purpose, let $\{e_1,\cdots,e_n\}\subset T_{y_0}M$ be an orthonormal basis for $T_{y_0}M$, i.e. it is a basis for $T_{y_0}M$ with $\langle e_i,e_j\rangle=\delta_{ij}$ for $i,j=1,\cdots, j$, where $\delta_{ij}$ is the Kronecker delta symbol. Denote by $e_i(t)=L_{y_0\bar y(t)}^{\bar y(\cdot)}e_i$ for $i=1,\cdots,n$ and $t\in(0,T]$, which is the parallel translation of $e_i$  from $y_0$ to $\bar y(t)$ along the curve $\bar y(\cdot)$ (See \cite[Section 2.2]{cdz}). We obtain from \cite[(2.6)]{cdz} that $\{e_1(t),\cdots,e_n(t)\}$ forms an orthonormal basis for $T_{\bar y(t)}M$, i.e. $\langle e_i(t),e_j(t)\rangle=\delta_{ij}$ for all $i,j\in\{1,\cdots,n\}$. Thus, we can write
\begin{align}\label{epf}
L_{y_k(t)\bar y(t)}f(t,y_k(t),u_k(t))=\sum_{i=1}^nf_k^i(t)e_i(t),\quad\forall\,k\geq\bar k\;\textrm{and}\;a.e. \,t\in[0,T],
\end{align}
where $f_k^i(t)\stackrel{\triangle}{=}\langle L_{y_k(t)\bar y(t)}f(t,y_k(t),u_k(t)),e_i(t)\rangle$. Set $\mathbb F_k(t)=(f_k^1(t),\cdots,f_k^n(t))^\top$, a.e. $t\in[0,T]$. It follows from \cite[(2.6)]{cdz} that
\begin{align}\label{ef1}
|f(t,y_k(t),u_k(t))|^2=|L_{y_k(t)\bar y(t)}f(t,y_k(t),u_k(t))|^2=\left|\sum_{i=1}^n f_k^i(t)e_i(t)\right|^2=\sum_{i=1}^nf_k^i(t)^2
\end{align}
holds for $k\geq \bar k$ and  a.e. $t\in[0,T]$.

\medskip

{\it Step 3.}\;We now analyze the convergence of the sequence $\{L_{y_k(\cdot)\bar y(\cdot)}f(\cdot,y_k(\cdot), u_k(\cdot))\}$ in a suitable sense. By means of (\ref{ef1}) and (\ref{ef}), we derive that
\begin{align*}
\int_0^T |\mathbb F_k(t)|^pdt&=\int_0^T \big(\sum_{i=1}^nf_k^i(t)^2\big)^{p/2}dt
\\&=\int_0^T|f(t,y_k(t),u_k(t))|^pdt
\\&\leq\int_0^T\big(K\rho(y_k(t),x_0)+\ell(t)\big)^pdt.
\end{align*}
Recalling that $\{y_k(\cdot)\}$ is bounded in $C([0,T];M)$, we conclude that $\{\mathbb F_k(\cdot)\}_{k\geq \bar k}$ is bounded in $L^p(0,T;\mathbb R^n)$. Thus, there exists a subsequence of $\{\mathbb F_k(\cdot)\}_{k\geq \bar k}$ (we still denote it by $\{\mathbb F_k(\cdot)\}_{k\geq \bar k}$) and $\bar{\mathbb F}(\cdot)=( \bar f^1(\cdot),\cdots,\bar f^n(\cdot))^\top\in L^p(0,T;\mathbb R^n)$ such that $\{\mathbb F_k(\cdot)\}_{k\geq\bar k}$ converges weakly to $\bar{\mathbb F}(\cdot)$ in $L^p(0,T;\mathbb R^n)$. Applying Mazur's theorem to $\{\mathbb F_k(\cdot)\}_{k\geq\bar k}$, we obtain that, for each $k\geq\bar k$, there exist $\beta_{k1},\cdots,\beta_{k n_k}\in[0,1]$ for some $n_k\in\mathbb N$ such that $\sum_{i=1}^{n_k}\beta_{ki}=1$ and the sequence $\{\sum_{i=1}^{n_k}\beta_{ki}\mathbb F_{k+i}(\cdot)\}_{k\geq \bar k}$ converges to $\bar{\mathbb F}(\cdot)$ strongly in $L^p(0,T;\mathbb R^n)$. Hence there exists a subsequence of  $\{\sum_{i=1}^{n_k}\beta_{ki}\mathbb F_{k+i}(\cdot)\}_{k\geq \bar k}$ (we still denote it by $\{\sum_{i=1}^{n_k}\beta_{ki}\mathbb F_{k+i}(\cdot)\}_{k\geq \bar k}$)
 such that
 \begin{align*}
\lim_{k\to+\infty}\sum_{j=1}^{n_k}\beta_{kj}\mathbb F_{k+j}(t)=\bar{\mathbb F}(t),\quad \hbox{a.e.}\;t\in[0,T].
\end{align*}
Set
\begin{align*}
\bar f(t)=\sum_{i=1}^n\bar f^i(t)e_i(t),\quad\bar f^0(t)={\lim\inf}_{k\to+\infty}\sum_{j=1}^{n_k}\beta_{kj}f^0(t,y_{k+j}(t),u_{k+j}(t)),\quad \hbox{a.e.}\;t\in[0,T].
\end{align*}
  Then, for a.e. $t\in[0,T]$,
\begin{align*}
\bar f(t)&=\sum_{i=1}^n\lim_{k\to+\infty}\sum_{j=1}^{n_k}\beta_{kj}f_{k+j}^i(t)e_i(t)
\\&=\lim_{k\to+\infty}\sum_{j=1}^{n_k}\beta_{kj}\sum_{i=1}^n f_{k+j}^i(t)e_i(t)
\\&=\lim_{k\to+\infty}\sum_{j=1}^{n_k}\beta_{kj}L_{y_{k+j}(t)\bar y(t)}f(t,y_{k+j}(t),u_{k+j}(t)).
\end{align*}
 Fix any $\delta\in(0,\bar \delta)$. By the convergence of $\{y_k(\cdot)\}$, we obtain that, there exists $k_\delta\geq \bar k$ such that
$
\rho(\bar y(t),y_k(t))<\delta$ holds for  every $t\in[0,T]$ and $k\geq k_\delta$.  Thus,
\begin{align*}
(\bar f^0(t),\bar f(t))\in \textrm{cl\,co}\bigcup_{\rho(z,\bar y(t))<\delta}L_{z\bar y(t)}\mathcal E(t,z),\quad \hbox{a.e.}\;t\in[0,T].
\end{align*}
Since $\delta\in(0,\bar \delta)$ is arbitrarily chosen, by the assumption (A6), we obtain that
\begin{align}\label{eff0}
\ds(\bar f^0(t),\bar f(t))\in\bigcap_{0<\delta<i(\bar y(t))} \textrm{cl\,co}\bigcup_{\rho(z,\bar y(t))<\delta}L_{z\bar y(t)}\mathcal E(t,z)=\mathcal E(t,\bar y(t)),\quad \hbox{a.e.}\;t\in[0,T].
\end{align}

\medskip

{\it Step 4.}\; In the step, we shall find an admissible control $\bar u(\cdot)$ such that $(\bar u(\cdot),\bar y(\cdot))$ is an admissible pair.

First, we  apply \cite[Theorem 8.2.9]{AF} to maps $(t,u)\mapsto f(t,\bar y(t),u)$, $t\mapsto \Gamma(t,\bar y(t))$ and $t\mapsto \bar f(t)$.  For every $u\in U$, we obtain from the assumption (A2) that the map $(t,y)\mapsto f(t,y,u)$ satisfies
\begin{description}
\item[(i)] For a.e. $t\in[0,T]$, the map $y(\in M)\mapsto f(t,y,u)(\in T M)$ is continuous;
\item[(ii)] For every $y\in M$, the map $t(\in[0,T])\mapsto f(t,y,u)(\in T M)$ is measurable
\end{description}
Applying \cite[Lemma 8.2.3, p. 311]{AF} to the maps $(t,y)\mapsto f(t,y,u)$ and $t\mapsto \bar y(t)$, we obtain that the map $t\mapsto f(t,\bar y(t),u)$ is measurable. On the other hand, for a.e. $t\in[0,T]$, we obtain from (A2) that $u\mapsto f(t,\bar y(t),u)$ is continuous.

Applying \cite[Theorem 8.2.8, p. 314]{AF} to $\Gamma(\cdot,\cdot)$ and  $\bar y(\cdot)$, we obtain from (A4) that, the map $t\mapsto \Gamma(t,\bar y(t))$ is measurable. By means of \cite[Theorem 8.2.9, p. 315]{AF}, there exists $\bar u(\cdot)\in\mathcal U$ such that
 \begin{align}\label{efc}
\bar u(t)\in\Gamma(t,\bar y(t))\quad\textrm{and}\quad \bar f(t)=f(t,\bar y(t),\bar u(t)),\quad \hbox{a.e.}\;t\in[0,T].
\end{align}
Consequently,
\begin{align}\label{ef01}
\bar f^0(t)\geq f^0(t,\bar y(t),\bar u(t)),\quad \hbox{a.e.}\;t\in[0,T].
\end{align}

Next, we claim that
\begin{align}\label{230701e1}
\dot{\bar y}(t)=f(t,\bar y(t),\bar u(t)),\quad \hbox{a.e.}\;t\in[0,T],
\end{align}
which immediately implies $(\bar u(\cdot),\bar y(\cdot))\in\mathcal P_{ad}$.
In fact,  by the compactness of $\big\{\bar y(t)\;\big|\;t\in[0,T]\big\}$, there exist $N\in\mathbb N$ and $t_0,t_1,\cdots,t_N \in [0,T]$ with $0=t_0<t_1<\cdots<t_N=T$ such that
\begin{align*}
\rho(\bar y(t_i),\bar y(t))<\bar\delta/2,\quad\forall\,t\in[t_i,t_{i+1}]\;\textrm{with}\;i\in\{0,1,\cdots,N-1\}.
\end{align*}
Fix any $i\in \{0,1,\cdots,N-1\}$.  By the above relation we can define $\exp_{\bar y(t_i)}^{-1}\bar y(t)$ and $\exp_{\bar y(t_i)}^{-1}y_k(t)$ for $t\in[t_i,t_{i+1}]$, when $k$ is large enough. Take any $X\in T_{\bar y(t_i)}M$.  For any $t\in[t_i,t_{i+1}]$, by means of \cite[Proposition 5.5.1, p. 187]{pth}, \cite[(5.13) and (5.15)]{dzgJ} and integration by parts over  $[t_i,t]$, we derive that
\begin{align}\label{ped4494}
\begin{array}{ll}
&\langle\exp_{\bar y(t_i)}^{-1}\bar y(t),X\rangle
\\[2mm]&=\ds\lim_{k\to+\infty}\langle\exp_{\bar y(t_i)}^{-1}y_k(t),X\rangle
\\&=\ds-\frac{1}{2}\lim_{k\to+\infty}\Big(\nabla_1\rho^2(\bar y(t_i),y_k(t))(X)-\nabla_1\rho^2(\bar y(t_i),y_k(t_i))(X) \Big)
\\[2mm]&=\ds-\frac{1}{2}\lim_{k\to+\infty}\int_{t_i}^t\nabla_2\nabla_1\rho^2(\bar y(t_i),y_k(s))\big(X, f(s,y_k(s),u_k(s))\big)ds
\\[2mm]&\ds=-\frac{1}{2}\int_{t_i}^t\nabla_2\nabla_1\rho^2(\bar y(t_i),\bar y(s))\big(X,\bar f(s)\big)ds+r_1(t)+r_2(t)
\\[2mm]&=\ds\Big\langle\int_{t_i}^t d\exp_{\bar y(t_i)}^{-1}|_{\bar y(s)}\bar f(s)ds, X\Big\rangle+r_1(t)+r_2(t),
\end{array}
\end{align}
where
\begin{align*}
r_1(t)\stackrel{\triangle}{=}&-\frac{1}{2}\lim_{k\to+\infty}\int_{t_i}^t\nabla_2\nabla_1\rho^2(\bar y(t_i),\bar y(s))\big(X,L_{y_k(s)\bar y(s)} f(s,y_k(s),u_k(s))-\bar f(s)\big)ds,
\\ r_2(t)\stackrel{\triangle}{=}&-\frac{1}{2}\lim_{k\to+\infty}\int_{t_i}^t\Big(\nabla_2\nabla_1\rho^2(\bar y(t_i),y_k(s))\big(X, f(s,y_k(s),u_k(s))\big)
\\ &-\nabla_2\nabla_1\rho^2(\bar y(t_i),\bar y(s))\big(X,L_{y_k(s)\bar y(s)}f(s,y_k(s),u_k(s))\Big)ds,
\end{align*}
and the integral in the last line is the Lebesgue integral on the Euclidean  space $T_{\bar y(t_i)}M$ with metric $\langle\cdot,\cdot\rangle$ evaluated at $\bar y(t_i)$. By (\ref{epf}), the definition of $\bar f(\cdot)$, and the weak convergence of $\{\mathbb F_k(\cdot)\}$ in $L^p(0,T;\mathbb R^n)$, we derive that
 \begin{align*}
r_1(t)=-\frac{1}{2}\lim_{k\to+\infty}\int_{t_i}^t\sum_{j=1}^n(f_k^j(s)-\bar f^j(t))\nabla_2\nabla_1\rho^2(\bar y(t_i),\bar y(s))(X,e_j(s))ds=0.
\end{align*}
For the term $r_2(t)$, we obtain from  \cite[(5.5) and (5.15)]{dzgJ}, the convergence of $\{y_k(\cdot)\}$ in the space $C([0,T];M)$, (\ref{ef1}), the boundedness of $\{\mathbb F_k(\cdot)\}$ in $L^p(0,T;\mathbb R^n)$, and \cite[Lemma 5.1]{dzgJ} that
\begin{align*}
r_2(t)&=\lim_{k\to+\infty}\int_{t_i}^t\Big(\langle d\exp_{y_k(s)}^{-1}|_{\bar y(t_i)}X,f(s,y_k(s),u_k(s))\rangle
\\&\quad-\langle d\exp_{\bar y(s)}^{-1}|_{\bar y(t_i)}X, L_{y_k(s)\bar y(s)}f(s,y_k(s),u_k(s))\rangle\Big)ds
\\&=\lim_{k\to+\infty}\int_{t_i}^t\langle d\exp_{y_k(s)}^{-1}|_{\bar y(t_i)}X-L_{\bar y(s)y_k(s)} d\exp_{\bar y(s)}^{-1}|_{\bar y(t_i)}X,f(s,y_k(s),u_k(s))\rangle ds
\\&=0.
\end{align*}
Note that $X\in  T_{\bar y(t_i)}M$ is arbitrarily chosen, we obtain that
\begin{align*}
\exp_{\bar y(t_i)}^{-1}\bar y(t)=\int_{t_i}^t d\exp_{\bar y(t_i)}^{-1}|_{\bar y(s)}\bar f(s)
ds,\quad\forall\,t\in[t_i,t_{i+1}],
\end{align*}
which is a curve in $T_{\bar y(t_i)}M$. This  curve is differentiable for a.e. $t\in[t_i,t_{i+1}]$. Thus, the curve
\begin{align*}
\bar y(t)=\exp_{\bar y(t_i)}\int_{t_i}^t d\exp_{\bar y(t_i)}^{-1}|_{\bar y(s)}\bar f(s)ds,\quad\forall \,t\in[t_i,t_{i+1}]
\end{align*}
is differentiable  for a.e. $t$, and
\begin{align*}
\dot{\bar y}(t)= d\exp_{\bar y(t_i)}|_{\exp_{\bar y(t_i)}^{-1}\bar y(t)}\circ d\exp_{\bar y(t_i)}^{-1}|_{\bar y(t)}\bar f(t),\quad \hbox{a.e.}\;t\in[t_i,t_{i+1}].
\end{align*}
Applying \cite[Proposition 3.6, p. 55]{le1} to $\exp_{\bar y(t_i)}$ and $\exp_{\bar y(t_i)}^{-1}$, we obtain that $\dot{\bar y}(t)=\bar f(t)$ for a.e. $t\in[t_i,t_{i+1}]$, which  together with (\ref{efc}), confirms (\ref{230701e1}).

Further,  we obtain from the assumption (A5) and the convergence of $\{y_k(\cdot)\}$ in $C([0,T];$ $M)$ that $(\bar y(0),\bar y(T))\in S$ and $\bar y(t)\in Q$ for all $t\in[0,T]$. Thus, $(\bar u(\cdot),\bar y(\cdot))\in\mathcal P_{ad}$.

\medskip

{\it Step 5.}\; By means of the assumption (A3) and (\ref{ef01}), one has
\begin{align*}
J(\bar u(\cdot),\bar y(\cdot))&=h(\bar y(0),\bar y(T))+\int_0^T f^0(t,\bar y(t),\bar u(t))dt
\\ &\ds\leq\liminf_{k\to+\infty}h( y_k(0),y_k(T))+\int_0^T\bar f^0(t)dt
\\&=\liminf_{k\to+\infty}h( y_k(0),y_k(T))+\int_0^T\liminf_{k\to+\infty}\sum_{j=1}^{n_k}\beta_{kj}f^0(t,y_{k+j}(t),u_{k+j}(t))dt
\\&\leq\liminf_{k\to+\infty}h( y_k(0),y_k(T))+\limsup_{k\to+\infty}\int_0^Tf^0(t,y_k(t),u_k(t))dt
\\&\leq\limsup_{k\to+\infty}\Big(h(y_k(0),y_k(t))+\int_0^Tf^0(t,y_k(t),u_k(t))dt\Big)
\\&=\lim_{k\to+\infty}J(u_k(\cdot),y_k(\cdot)),
\end{align*}
which implies $(\bar u(\cdot),\bar y(\cdot))$ is an optimal pair.   The proof of Theorem \ref{eth} is completed.
$\Box$

\medskip

\textbf{Proof of Proposition \ref{vctc}}\; It suffices to prove that, $\mathcal E(t,\cdot)$ satisfies the Cesari-type property on $M$ provided that $\mathcal E(t,x)$ is convex and closed for each $x\in M$.

 For the case that (A$3^\prime$) holds, we fix any $\epsilon>0$. It follows from the condition
 (A$3^\prime$)  that, there exists $\sigma\in(0,i(x))$ such that, for any $(x^\prime,u,u^\prime)\in M\times U\times U$ with $\max\{\rho(x,x^\prime),d(u,u^\prime)\}<\sigma$, one has
 \begin{align}\label{ctp11}
|L_{x^\prime x}f(t,x^\prime,u^\prime)-f(t,x,u)|<\frac{\epsilon}{\sqrt 2}\quad\textrm{and}\quad f^0(t,x^\prime,u^\prime)>f^0(t,x,u)-\frac{\epsilon}{\sqrt{2}}.
\end{align}
Since $\Gamma(t,\cdot)$ is upper semicontinuous, it follows from \cite[Proposition 2.2, p. 89]{ly} that, there exists $\delta^\prime\in(0,\sigma)$ such that
\begin{align*}
\Gamma(t,B_x(\delta^\prime))\subseteq B_\sigma(\Gamma(t,x)),
\end{align*}
where $B_x(\delta^\prime)\stackrel{\triangle}{=}\big\{z\in M\;\big|\;\rho(z,x)<\delta^\prime\big\}$ and  $B_\sigma(\Gamma(t,x))\stackrel{\triangle}{=}\big\{v\in U\;\big|\;\inf_{w\in\Gamma(t,x)}d(w,v)<\sigma\big\}$.
For any $z\in B_x(\delta^\prime)$, take any $(\zeta^0,\zeta)\in\mathcal E(t,z)$. Choose $u_z\in \Gamma(t,z)$ and $u^\prime_z\in \Gamma(t,x)$ such that
\begin{align*}
\zeta=f(t,z,u_z),\quad\zeta^0\geq f^0(t,z,u_z)\quad \textrm{and}\quad d(u_z,u_z^\prime)<\sigma.
\end{align*}
It follows from (\ref{ctp11}) that $(\zeta^0,L_{z x}\zeta)\in B_\epsilon(\mathcal E(t,x))$, where $B_\epsilon(\mathcal E(t,x))\stackrel{\triangle}{=}\big\{(y^0,Y)\in \mathbb R\times T_xM\;\big|\;\inf_{(z^0,Z)\in \mathcal E(t,x)}|z^0-y^0|^2+|Z-Y|^2<\epsilon^2\big\}$.
Since  both $z\in B_x(\delta^\prime)$ and $(\zeta^0,\zeta)\in \mathcal E(t,z)$ are arbitrarily chosen,  one has $\bigcup_{\rho(y,x)<\delta^\prime}L_{y x}\mathcal E(t,y)\subseteq B_\epsilon(\mathcal E(t,x))$, which implies that
\begin{align*}
\bigcap_{i(x)>\delta>0}\textrm{cl}\,\textrm{co}\,\big(\bigcup_{\rho(y,x)<\delta}L_{yx}\mathcal E(t,y)\big)\subseteq  \textrm{cl}\,\textrm{co}\,B_\epsilon(\mathcal E(t,x))\subseteq \textrm{cl}\,B_\epsilon(\mathcal E(t,x)).
\end{align*}
Here one has used the relation $\textrm{co}\, B_\epsilon(\mathcal E(t,x))\subseteq B_\epsilon(\textrm{co}\, \mathcal E(t,x))$ and the assumption that $\mathcal E(t,x)$ is convex. Since $\epsilon>0$ is arbitrarily chosen and $\mathcal E(t,x)$ is closed, (\ref{ctpx}) follows.

For the case that (A$4^\prime$) holds, the proof is analogous to the previous case.  This completes the proof of Proposition \ref{vctc}. $\Box$

\bibliographystyle{amsplain}
\providecommand{\bysame}{\leavevmode\hbox to3em{\hrulefill}\thinspace}
\providecommand{\MR}{\relax\ifhmode\unskip\space\fi MR }
\providecommand{\MRhref}[2]{%
  \href{http://www.ams.org/mathscinet-getitem?mr=#1}{#2}
}
\providecommand{\href}[2]{#2}

\end{document}